\documentclass{article}
\usepackage{fullpage}
\usepackage{graphics,psfrag,epsfig}
\usepackage{amsfonts,latexsym,eucal,amsmath,amsthm,amssymb}

\begin{document}

\newcommand{\rum}{\rule{0.5pt}{0pt}}
\newcommand{\rub}{\rule{1pt}{0pt}}
\newcommand{\rim}{\rule{0.3pt}{0pt}}
\newcommand{\numtimes}{\mbox{\raisebox{1.5pt}{${\scriptscriptstyle \times}$}}}
\newcommand{\optprog}[2]
{%
  \noindent\mbox{}\\[0cm]
  \noindent\fbox{%
  \begin{minipage}{0.955\linewidth}
    \mbox{}\\[-0.5cm]
    #1\\[#2]
  \end{minipage}
  }
  \noindent\mbox{}\\[-0.2cm]
}

\renewcommand{\refname}{References}

\twocolumn[%
\begin{center}
{\Large\bf The Riemann Hypothesis is Unprovable \rule{0pt}{13pt}}\par
\bigskip
Craig Alan Feinstein \\ {\small\it 2712 Willow Glen Drive, Baltimore, Maryland
21209\rule{0pt}{13pt}}\\ \raisebox{-1pt}{\footnotesize E-mail: cafeinst@msn.com,
BS"D}\par
\bigskip\smallskip
{\small\parbox{11cm}{%
\bigskip \noindent \textbf{Abstract:} In this note, we give a simple proof that the Riemann Hypothesis
is unprovable in any reasonable axiom system.

\bigskip \noindent \textbf{Disclaimer:} This article was authored
by Craig Alan Feinstein in his private capacity. No official support or endorsement by
the U.S. Government is intended or should be inferred.\rule[0pt]{0pt}{0pt}}}\bigskip
\end{center}]{%

The Riemann-Zeta function $\zeta(s)$ is a complex function defined to be
$$
\zeta(s)=\frac{s}{s-1}-s \int_1^{\infty}\frac{x-\lfloor x \rfloor}{x^{s+1}} dx
$$
when the real part of the complex number $s$ is positive \cite{b:RZ}. The Riemann Hypothesis states
that if $\rho=\sigma+ti$ is a complex root of $\zeta(s)$ and
$0<\sigma<1$, then $\sigma=1/2$ \cite{b:RH}. The Riemann-Siegel function is defined to
be the real function
$$
Z(t)=\zeta(\frac{1}{2} + ti)\cdot \exp(i\vartheta(t)),
$$
where
$$
\vartheta(t)= \arg [\Gamma ( \frac{1}{4} + \frac{1}{2}it) ] - \frac{1}{2} t \ln \pi.
$$
Notice that $|Z(t)|=|\zeta(1/2 + ti)|$ for all real $t$, so the real roots
$t$ of $Z(t)$ are the same as the real roots $t$ of $\zeta(1/2 + ti)$ \cite{b:RS}.
In this note, we give a simple proof that the Riemann Hypothesis is unprovable in
any reasonable axiom system:

The Riemann Hypothesis is equivalent to the assertion that for each $T>0$, the number of real roots $t$ of $\zeta(1/2+it)$ (counting multiplicities), where $0<t<T$, is equal to the number of complex roots $s$ of $\zeta(s)$ in $\{s=\sigma+ti \mid 0<\sigma<1,\mbox{ } 0<t<T\}$ (counting multiplicities). Because the formula for the real roots $t$ of $\zeta(1/2+it)$ cannot be reduced to a formula that is simpler than the equation, $\zeta(1/2+it)=0$, the only way to determine the number of real roots $t$ of $\zeta(1/2+it)$ in which $0<t<T$ is to count the changes in sign of $Z(t)$, where $0<t<T$.

Suppose it were possible to prove that the number of real roots $t$ of $\zeta(1/2+it)$,
where $0<t<T$, is equal to the number of complex roots $s$ of $\zeta(s)$ in $\{s=\sigma+ti\mid 0<\sigma<1,\mbox{ } 0<t<T\}$ without counting the changes in sign of $Z(t)$, where $0<t<T$. Then it
would be possible to determine the number of real roots $t$ of $\zeta(1/2+it)$, where $0<t<T$, by computing the number of complex roots $s$ of $\zeta(s)$ in $\{s=\sigma+ti\mid 0<\sigma<1,\mbox{ } 0<t<T\}$ via the argument principle \cite{b:AP} without counting the changes in sign of $Z(t)$, where $0<t<T$. But this contradicts the last sentence in the paragraph above. Hence, in order to prove that the number of real roots $t$ of $\zeta(1/2+it)$, where $0<t<T$, is equal to the number of complex roots $s$ of $\zeta(s)$ in $\{s=\sigma+ti\mid 0<\sigma<1,\mbox{ } 0<t<T\}$, it is necessary to count the changes in sign of $Z(t)$, where $0<t<T$.

As $T$ becomes arbitrarily large, the time that it takes to count the changes in sign of $Z(t)$, where $0<t<T$, approaches infinity; hence, an infinite amount of time is required to prove that for each $T>0$, the number of real roots $t$ of $\zeta(1/2+it)$, where $0<t<T$, is equal to the number of complex roots $s$ of $\zeta(s)$ in $\{s=\sigma+ti \mid 0<\sigma<1,\mbox{ } 0<t<T\}$, so the Riemann Hypothesis is unprovable in any reasonable axiom system.

}

\end{document}